\documentclass[UTF-8,reqno,12pt]{article}
\usepackage{amsmath}
\usepackage{amssymb}
\usepackage{amsthm}
\usepackage[libertinus]{newtx}  
\usepackage{geometry}
\usepackage{xcolor}

\usepackage[colorlinks=true,backref=page]{hyperref}
\hypersetup{
    linkcolor=blue,          
    citecolor=red,        
    filecolor=blue,      
    urlcolor=magenta   
}
\usepackage[nobysame,abbrev]{amsrefs}
\renewcommand{\eprint}[1]{{\it Available at}\href{https://arxiv.org/abs/#1}{~\color{blue}\it{arXiv:#1}}.}
\renewcommand{\PrintDOI}[1]{\url{https://doi.org/#1}}
\renewcommand{\MR}[1]{\href{https://mathscinet.ams.org/mathscinet-getitem?mr=#1}{\color{cyan}{MR#1}}}

\BibSpec{article}{%
+{}{\PrintAuthors} {author}
+{,}{ \textrm} {title}
+{:}{ \textrm} {subtitle}
+{.}{ \textit} {journal}
+{,}{ \textbf} {volume}
+{}{ \parenthesize} {date}
+{,}{ } {pages}
+{}{ \parenthesize} {language}
+{.}{} {transition}
+{.}{ \PrintReviews} {review}
+{.}{ \eprint} {eprint}
+{.}{ \PrintDOI} {doi}
+{.}{ \url} {url}
}
\BibSpec{book}{%
+{}{\PrintAuthors} {author}
+{,}{ \textit} {title}
+{:}{ \textit} {subtitle}
+{.}{ \textrm} {series} 
+{,}{ Vol.} {volume} 
+{.}{ } {publisher}
+{,}{ } {date}
+{}{ \parenthesize} {language}
+{.}{} {transition}
+{.}{ \PrintReviews} {review}
+{.}{ \PrintDOI} {doi}
}
\theoremstyle{plain}
\newtheorem{theorem}{Theorem}[section]

\newtheorem{proposition}[theorem]{Proposition}
\newtheorem{lemma}[theorem]{Lemma}
\theoremstyle{definition}
\newtheorem{definition}[theorem]{Definition} 
\numberwithin{equation}{section}
\theoremstyle{remark}
\newtheorem{remark}[theorem]{Remark} 
\newtheorem{example}[theorem]{Example}
\theoremstyle{Assumption}
\def\bt{\begin{theorem}}
\def\et{\end{theorem}}
\def\bl{\begin{lemma}}
\def\el{\end{lemma}}
\def\bd{\begin{definition}}
\def\ed{\end{definition}}
\def\bc{\begin{corollary}}
\def\ec{\end{corollary}}
\def\br{\begin{remark}}
\def\er{\end{remark}}
\def\bp{\begin{proposition}}
\def\ep{\end{proposition}}
\def\bexa{\begin{example}}
\def\eexa{\end{example}}

\def\bB{{\mathbf B}}
\def\bC{{\mathbf C}}
\def\bH{{\mathbf H}}
\def\bP{{\mathbf P}}
\newcommand\1{{\boldsymbol{1}}}

\def\B{\mathrm{B}}

\def\cP{{\mathcal P}}
\def\cR{{\mathcal R}}

\def\mE{{\mathbb E}}
\def\mI{{\mathbb I}}
\def\mM{{\mathbb M}}
\def\mN{{\mathbb N}}
\def\mP{{\mathbb P}}
\def\mR{{\mathbb R}}
\def\mS{{\mathbb S}}

\def\sA{{\mathscr A}}
\def\sB{{\mathscr B}}
\def\sI{{\mathscr I}}
\def\sL{{\mathscr L}}
\def\sF{{\mathscr F}}
\def\sS{{\mathscr S}}

\def\p{\partial}
\def\geq{\geqslant}
\def\leq{\leqslant}
\def\ge{\geqslant}
\def\le{\leqslant}
\def\eps{\varepsilon}

\def\[{{\Big[}}
\def\]{{\Big]}}
\def\<{{\langle}}
\def\>{{\rangle}}
\def\({{\Big(}}
\def\){{\Big)}}

\def\e{\mathrm{e}}
\providecommand{\dif}{\mathop{}\!\mathrm{d}} 

\def\div{{\mathord{\mathrm{div}}}}

\def\bpf{\begin{proof}}
\def\epf{\end{proof}}

\allowdisplaybreaks  
\sloppypar

\date{ }

\begin{document}

\title{Euler--Maruyama scheme for $\alpha$-stable SDE with distributional drift}


\author{Zimo Hao \thanks{School of Mathematics and Statistics, Beijing Institute of Technology, Beijing, 100081,  P. R. China. Email: zimo\_hao@163.com}~~~and~~
Mingyan Wu \thanks{School of Mathematical Sciences, Xiamen University, Xiamen, Fujian, 361005, P. R. China. E-mail: mingyanwu.math@xmu.edu.cn; mingyanwu.math@gmail.com}}

\maketitle

\begin{abstract}
In this paper, we consider a class of stochastic differential equations driven by symmetric non-degenerate $\alpha$-stable processes (including cylindrical ones) with $\alpha \in (1,2)$. We first establish a quantitative estimate for the Euler scheme under bounded drift $b(x)$, with an explicit dependence on $ \| b \|_{L^\infty}$. Then we obtain the weak convergence rates for the case where the drift coefficient belongs to a Besov space of negative order.
\end{abstract}



\tableofcontents

\section{Introduction}

Recently, stochastic differential equations (SDEs) with distributional drifts have attracted considerable attention, both for Brownian noise (see e.g., \cite{DD16,CC18,HZ23}) and for $\alpha$-stable noise (see e.g., \cite{ABM18,CM19,KP20,LZ22}). Beyond motivations arising from regularization by noise, SDEs with distributional drifts often model random irregular media and exhibit distinct behaviors. Examples include Brox diffusion (see \cite{HLM17}), superdiffusive phenomena \cite{CHT22,CMOW23}, random directed polymers \cite{DD16}, and self-attracting Brownian motion in a random medium \cite{CC18}. For further references on the motivations for studying SDEs with distributional drifts, we refer the reader to \cite{DGI22}.

\medskip
 In this paper, we investigate the Euler--Maruyama approximation of the following SDE in $\mR^d$ ($d\geq 1$):
\begin{align}\label{in:SDE}
\dif X_t = b(X_t) \dif t + \dif L_t^{(\alpha)},\qquad X_0 = x \in \mR^d,
\end{align}
where the drift coefficient $b$ belongs to $\bB_{\infty,\infty}^{-\beta}(\mR^d)$ for some $\beta \in (0, \alpha-1)$ (here, $\bB_{\infty,\infty}^{-\beta}$ denotes a Besov space; see Definition \ref{iBesov} below), and $L^{(\alpha)}$ is a $d$-dimensional symmetric $\alpha$-stable process with $\alpha \in (1,2)$ on some probability space $(\Omega, \sF, \mP)$. Its L\'evy measure is given by
\begin{align}\label{CC02}
\nu^{(\alpha)}(A) = \int_0^\infty \left( \int_{\mS^{d-1}} \frac{1_A(r\theta) \, \Sigma(\dif\theta)}{r^{1+\alpha}} \right) \dif r,\qquad A \in \sB(\mR^d),
\end{align}
where $\Sigma$ is a finite measure on the unit sphere $\mS^{d-1}$. This formulation unifies two important cases:
\begin{itemize}
\item If $\Sigma$ is the uniform (rotation-invariant) measure on $\mS^{d-1}$, then $L^{(\alpha)}$ is the \emph{standard} (rotationally invariant) $\alpha$-stable process. Its L\'evy measure is absolutely continuous with respect to the Lebesgue measure, given by $\frac{1}{|z|^{d+\alpha}} \dif z$, and its infinitesimal generator is the fractional Laplace operator $\Delta^{\alpha/2}$. Notice that the components of a standard $\alpha$-stable process are not jointly independent.
\item If $\Sigma$ is concentrated on the coordinate axes, i.e., $\Sigma = \sum_{i=1}^d \delta_{\pm e_i}$, then $L^{(\alpha)}$ becomes a \emph{cylindrical} $\alpha$-stable process, whose components are independent one-dimensional $\alpha$-stable processes. In this case, the L\'evy measure is given by
$$
\nu^{(\alpha)}(\dif z) := \sum_{k=1}^d \delta_0(\dif z_1) \cdots \delta_0(\dif z_{k-1}) \, \frac{\dif z_k}{|z_k|^{1+\alpha}} \, \delta_0(\dif z_{k+1}) \cdots \delta_0(\dif z_d),
$$
where $\delta_0$ is the Dirac measure at zero. Consequently, the symbol of its infinitesimal generator is $\sum_{i=1}^d |\xi_i|^\alpha$, which is more singular than that of the standard process: while $|\xi|^\alpha$ is non-smooth only at the origin, $\sum_{i=1}^d |\xi_i|^\alpha$ fails to be smooth on the entire set of coordinate axes $\bigcup_{i=1}^d \{\xi_i = 0\}$. This is why the cylindrical process is referred to as singular.
\end{itemize}
We point out that the joint independence of the components $\{L^i\}_{i=1}^d$ plays a vital role in many models. For instance, in the following $N$-particle system:
$$
\dif X^{N,i}_t = \frac{1}{N} \sum_{j \neq i} K(X^{N,i}_t - X^{N,j}_t) \dif t + \dif L^i_t,
$$
$K: \mR^d \to \mR^d$ is the interaction kernel, and $\{L^i\}_{i=1}^N$ is a family of independent $\alpha$-stable processes, which models random phenomena such as collisions between two particles (see \cite{Ca22} and references therein).

\medskip

Compared to the function-drift case,  only a few works concern numerical schemes for SDEs with distributional drifts. To the best of our knowledge, only three works (see \cite{DGI22, GHR25, CIP25}) have studied Euler-type approximations within the distributional framework. Specifically, \cite{DGI22} and \cite{CIP25} investigate the numerical solution of one-dimensional SDEs with distributional drifts and Brownian noise. The former considers drifts in fractional Sobolev spaces of negative regularity, while the latter treats drifts in Besov spaces of negative order. Additionally, \cite{GHR25} studies a tamed Euler scheme for $d$-dimensional SDEs with drifts in negative Besov spaces and noise given by fractional Brownian motion. It is worth pointing out that all the aforementioned works only establish strong convergence rates for continuous noise. No results on convergence rates are currently available for the case of $\alpha$-stable noise, even for the standard ones.

\medskip

In this work, we aim to fill this gap by developing a unified framework for the Euler--Maruyama approximation of SDEs driven by a class of $\alpha$-stable processes that includes both standard and cylindrical cases, with distributional drifts. The detailed problem statement and our main results are presented in Sections \ref{Sec:pp} and \ref{sec:main}, respectively.

\subsection*{Conventions and notations} 
Throughout this paper, we use the following conventions and notations: As usual, we use $:=$ as a way of definition. Define $\mN_0:= \mN \cup \{0\}$ and $\mR_+:=[0,\infty)$. The letter $c=c(\cdots)$ denotes an unimportant constant, whose value may change in different places. We use $A \asymp B$ and $A\lesssim B$ to denote $c^{-1} B \leq A \leq c B$ and $A \leq cB$, respectively, for some unimportant constant $c \geq 1$. 
Denote the Beta function by 
\begin{align}\label{eq:Beta}
\B (s_1,s_2):= \int_0^1 x^{s_1-1}(1-x)^{s_2-1} \dif x, \ \ \forall s_1,s_2>0.
\end{align}

\begin{itemize}
\item  Let $\mM^d$ be the space of all real $d\times d$-matrices, and $\mM^d_{non}$ the set of all non-singular matrices. Denote the identity $d\times d$-matrix by $\mI_{d \times d}$.
 
\item  For every $p\in [1,\infty)$, we denote by $L^p$ the space of all $p$-order integrable functions on $\mR^d$ with the norm denoted by $\|\cdot\|_p$. 

\item The norm $\|\cdot\|_\infty$ is defined as $\|f\|_\infty := \mathrm{ess\,sup}_{x \in \mR^d} |f(x)|$.
\item  Let $\mathcal{P}(\mathbb{R}^d)$ denote the set of all probability measures on $\mathbb{R}^d$.
\item Let $\|\mu_1 - \mu_2\|_{\rm var}$ denote the total variation distance between two probability measures $\mu_1$ and $\mu_2$ on $\mR^d$, defined by
\begin{align*}
\|\mu_1 - \mu_2\|_{\rm var} := \sup_{\|\varphi\|_\infty = 1} \left| \int_{\mR^d} \varphi(x) \, (\mu_1 - \mu_2)(\dif x) \right|.
\end{align*}
\end{itemize}

\subsection*{Organization of the paper}

The remainder of this paper is organized as follows. Section \ref{Sec:pp} states the problem and explains the transition from smooth to distributional coefficients. Section \ref{sec:main} presents our two main results. Section \ref{Pre} collects preliminaries on Besov spaces, $\alpha$-stable processes, and heat kernel estimates. Section \ref{sec:proof} establishes the weak convergence rates of the Euler scheme, first for bounded drifts (see Theorem \ref{thmEB}) and then for distributional drifts (see Theorem \ref{in:Main}).


\section{Problem statement}
\label{Sec:pp}

Since the drift term $b$ is a distribution, which is not meaningful in the classical sense, it is impossible to assign a value to a distribution at the point $X_t$. To define solutions and their Euler's scheme, a natural approach is to use mollifying approximations. Let $\phi_m(x) := m^d \phi(mx)$, $m \in\mN$, be a family of mollifiers, where $\phi \in C_c^\infty(\mathbb{R}^d)$ is a smooth probability density function with compact support. The smooth approximation of $b$ is then defined by convolution as follows:  
\begin{align}\label{BN1}
b_m (x) :=  ( b  * \phi_m )(x).
\end{align}

\medskip

We consider a mollified Euler's scheme for SDE \eqref{in:SDE}. Let $X^{m}_t$ solve the classical SDE
\begin{align*}
X^{m}_t = x + \int_0^t b_m(X^{m}_s) \dif s + L_t^{(\alpha)},
\end{align*}  
and $X^{m,n}_t$ be its Euler scheme: for any $n\in\mN$,
\begin{align}\label{0608:00}
X^{m,n}_t=x+\int_0^tb_m(X^{m,n}_{\pi_n(s)})\dif s+L_t^{(\alpha)},\ \ x \in \mR^d, t \in (0,T],
\end{align} 
where  $n \in \mN$, and $\pi_n(t):=k/n$ for $t\in[k/n,(k+1)/n)$ with $k=0,1,2,...., \lfloor n T \rfloor$. 
Thanks to the stability estimates (see Lemma \ref{thmSt}), to prove our main result on weak convergence rates of Euler's scheme (see Theorem \ref{in:Main}), it suffices to establish a quantitative estimate (Theorem \ref{thmEB}) for the difference between $X^{m}_t$ and $X^{m,n}_t$, with an explicit dependence on $\|b_m\|_\infty$. The key technique used in this task is the so-called It\^o--Tanaka trick, which has been widely used in the literature to obtain quantitative estimates of the Euler approximation for both continuous and discontinuous drifts (see e.g., \cite{TT90, MP91, Hol22, FJM25, SH24}). This trick exploits the regularizing effect of the semigroup. 

\medskip

Let us  first briefly recall the It\^o--Tanaka trick. Consider the function $u(t,x) := \mathbb{E} \varphi(x + L_t^{(\alpha)})$, which solves the PDE
\begin{align*} 
\p_t u=\sL^{(\alpha)} u,\quad u(0)=\varphi \in C^\infty_b,
\end{align*}
where 
$$
\sL^{(\alpha)} f (x) :=\int_{\mR^d}\(f(x+ z)-f(x)-  z \cdot\nabla f(x)\)\nu^{(\alpha)}(\dif z) .
$$
Applying It\^o's formula to $s \mapsto u(t-s, X^m_s)$ and $s \mapsto u(t-s, X^{m,n}_s)$ respectively, we obtain
\begin{align*}
|\mE \varphi(X^{m,n}_t)-\mE \varphi(X^{m}_t)|\leq &\left|\mE\int_0^t \Big((b_m\cdot\nabla u(t-s))(X^{m,n}_{s})-(b_m\cdot\nabla u(t-s))(X^{m}_{s})\Big)\dif s\right|\\
+&\left|\mE\int_0^t \left(b_m(X^{m,n}_{\pi_n (s)})-b_m(X^{m,n}_s)\right)\cdot \nabla u(t-s,X^{m,n}_s)\dif s\right|\\
\le& \int_0^t\|b_m\cdot \nabla u(t-s)\|_{\infty}\|\mP\circ(X^{m,n}_s)^{-1}-\mP\circ(X^{m}_s)^{-1}\|_{\rm var}\dif s\\
&+ \|b_m\|_{C_b^1}\mE\int_0^t \|\nabla u(t-s)\|_\infty|X^{m,n}_{\pi_n (s)}-X^{m,n}_s|\dif s.
\end{align*}
Since $\|\nabla u(t)\|_\infty\lesssim t^{-1/\alpha}\|\varphi\|_\infty$ (see \eqref{DH09}), taking the supremum over $\|\varphi\|_\infty = 1$ yields, for $\alpha > 1$ and $t \in (0,T]$, 
\begin{align*}
\|\mP\circ(X^{m,n}_t)^{-1} & -\mP\circ(X^{m}_t)^{-1} \|_{\rm var}
\lesssim t^{\frac{\alpha-1}{\alpha}} \|b_m\|_{ C_b^1} (\|b_m\|_\infty n^{-1}+n^{-\frac1\alpha}) \\
& +  \|b_m\|_\infty\int_0^t(t-s)^{-\frac1\alpha}\|\mP\circ(X^{m,n}_s)^{-1}-\mP\circ(X^{m}_s)^{-1}\|_{\rm var}\dif s,
\end{align*}
which, by Gronwall's inequality of  Volterra's type (see \cite{We19}, Theorem 3.2, or \cite{Zh10}, Lemma 2.2), derives that there are two constants $c_0=c_0(\|b_m\|_\infty)>0$ and $c_1=c_1(\|b_m\|_{C_b^1})>0$ such that for any $t \in (0,T]$ and $n \in \mN$, 
\begin{align}\label{0607:01}
\|\mP\circ(X^{m,n}_t)^{-1}-\mP\circ(X^{m}_t)^{-1}\|_{\rm var}\le  c_1\e^{c_0}t^{\frac{\alpha-1}{\alpha}} n^{-\frac1\alpha}.
\end{align}

\medskip
With the estimate \eqref{0607:01} in hand, we now have a quantitative control for the case with smooth coefficients. Returning to the original distributional setting, however, two main questions arise when we try to apply this estimate. Recall that $b \in \bB^{-\beta}_{\infty,\infty}$, and the mollified drift $b_m$, defined by \eqref{BN1}, satisfies 
\begin{align}\label{eq:FT00}
\|b_m\|_\infty \lesssim m^\beta \|b\|_{\bB^{-\beta}_{\infty,\infty}} \quad \text{and} \quad \|b_m\|_{C_b^1} \lesssim m^{\beta+1} \|b\|_{\bB^{-\beta}_{\infty,\infty}}.
\end{align}
In this context, we are led to the following two issues.

\medskip\noindent
\textbf{(1) Regularity of the drift.} 

\medskip

In the estimate \eqref{0607:01}, the constant $c_1$ depends positively on $\|b_m\|_{C_b^1}$, which grows like $m^{\beta+1}$ by \eqref{eq:FT00}. To minimize the growth of the mollification parameter $m$, we would like the dependence in $c_1$ to be on $\|b_m\|_\infty$ instead, which grows only like $m^\beta$. This raises the following question: can we reduce the dependence on $\|b_m\|_{C_b^1}$ in $c_1$ to a dependence on $\|b_m\|_\infty$?

\medskip

For this question, an initial qualitative result was given by Gy\"ongy and Krylov \cite{GK}, who showed that $X^{m,n}$ converges in probability to $X^{m}$ when the noise is Brownian motion and the drift $b$ is merely bounded and measurable. However, a quantitative result concerning the dependence on $\|b_m\|_\infty$ appears to be absent in the literature.

\medskip\noindent
\textbf{(2) Exponential growth.} 

\medskip

The factor $\e^{c_0}$ in \eqref{0607:01} grows like $\exp\bigl\{ m^{\frac{\alpha \beta}{\alpha-1}} \bigr\}$ (cf. Theorem 3.2 of \cite{We19}) since \eqref{eq:FT00}. To counteract this growth, one might choose $m \sim (\ln n)^{\frac{\alpha-1}{\alpha \beta}}$, where $n$ is the discretization parameter. However, this choice is not satisfactory for the following reasons:
\begin{itemize}
\item The mollification parameter $m$ grows only logarithmically in $n$, so an extremely large $n$ is required to make $m$ sufficiently large to ensure accurate approximation of the distributional drift;
\item Combining this choice with the stability estimates (see Lemma \ref{thmSt}) leads to a convergence rates that is logarithmic in $n$ rather than polynomial, which is too slow for practical purposes. In practice, one needs a polynomial relation between $m$ and $n$, e.g., $m = n^\gamma$ with $\gamma > 0$, to achieve a reasonable convergence rates.
\end{itemize}
This leads to the second question: can we obtain a polynomial dependence on $m$ instead of the exponential factor $\e^{c_0}$ in \eqref{0607:01}?

\medskip

To fix these two issues, we apply the It\^o--Tanaka trick twice. This allows us to obtain the desired estimates without relying on Gronwall's inequality, thereby avoiding both the dependence on $\|b_m\|_{C_b^1}$ and the exponential growth of the mollification parameter $m$. Consequently, we derive an upper bound that is polynomial in $n$ and depends explicitly on $\|b_m\|_\infty$ (see Theorem \ref{thmEB}). This leads to our second main result: the weak convergence rates for the Euler scheme of SDE \eqref{in:SDE} (see Theorem \ref{in:Main}) under the assumption $m = n^\gamma$ for some $\gamma > 0$.

\section{Main results}\label{sec:main}

Throughout this paper, we always assume that the following condition holds:

\medskip
\noindent
$\bf (ND)$ The L\'evy measure given by \eqref{CC02} is {\it non-degenerate}, that is, for each $ \theta_0\in \mS^{d-1}$,  $$
\int_{\mS^{d-1}}|\theta\cdot\theta_0|\Sigma(\dif \theta)>0.
$$

\br
Here, we refer to  \cite{HWW20}, Examples 2.10 and 2.11, as two examples of L\'evy processes satisfying the non-degeneracy  condition {\bf(ND)}. 
\er

We state the following definition of weak solutions to SDE \eqref{in:SDE}.

\begin{definition}[Weak solutions]\label{Def1}  
Let $ (\Omega, \sF, (\sF_t)_{t \geq 0}, \mathbb{P})$ be a stochastic basis, and let $(X, L)$ be a pair of $\mathbb{R}^d$-valued, c$\rm\grave{a}$dl$\rm\grave{a}$g, $(\sF_t)$-adapted processes on $ (\Omega, \sF, (\sF_t)_{t \geq 0}, \mathbb{P})$. We call $( X, L)$ with $ (\Omega, \sF, (\sF_t)_{t \geq 0}, \mathbb{P})$ a weak solution of the SDE \eqref{in:SDE} with initial distribution $\mu \in \mathcal{P}(\mathbb{R}^d)$ if $L$ is an $(\sF_t)$-$\alpha$-stable process with the L\'evy measure $\nu$ given by \eqref{CC02} which satisfies the condition $\bf (ND)$, and $\mathbb{P} \circ X_0^{-1} = \mu$, and  
$$
X_t = X_0 + A^b_t +   L_t,\quad \text{for all $t \in [0, T]$,} \quad \text{a.s.},
$$
where $A^b_t := \lim_{m \to \infty} \int_0^t b_m(X_s) \, \mathrm{d}s$ exists in the $L^2$-sense, with $b_m$ defined by \eqref{BN1}. 
\end{definition}


Fortunately, the well-posedness has been established by our previous work \cite{HW23}. For the reader's convenience, we present the result here.

\bp[Weak well-posedness]\label{thm:G-mart}
Let $T>0$, $\alpha\in(1,2)$ and $\beta\in(0,\alpha-1)$.
Assume that  
\begin{align}\label{con:b}
    (i)~~b\in \bB_{\infty,\infty}^{-\beta},\quad \text{if $\beta\in (0,\tfrac{\alpha-1}{2}$)}; \quad  \quad  (ii)~~b,\div b\in  \bB_{\infty,\infty}^{-\beta}, \quad \text{if $\beta\in [\tfrac{\alpha-1}{2}, \alpha-1)$}.
\end{align}
Then for any $\mu\in \cP(\mR^d)$, there is a unique weak solution to SDE \eqref{in:SDE} in the sense of Definition \ref{Def1}. The weak solution is independent of the specific choice of mollifier functions $\phi_m$.
\ep

For simplicity of notation,  we introduce the following parameter set:
\begin{align*}
\Theta:=(T,d,\alpha,\beta).
\end{align*}

\subsection{Quantitative estimates for bounded drift}

We first study the Euler scheme for SDEs with bounded drift. Let $n \in \mN$, $X_0^n = X_0 = x$, and define
\begin{align}\label{eq:JN01-0}
X^n_t = x + \int_0^t b(X^n_{\pi_n(s)}) \, \dif s + L_t^{(\alpha)}, \qquad t \in [0,T],
\end{align}
where $\pi_n(t) := k/n$ for $t \in [k/n, (k+1)/n)$ with $k = 0,1,2,\dots, \lfloor nT \rfloor$. Our first goal is to establish a quantitative estimate for the Euler scheme \eqref{eq:JN01-0} under bounded drift, where the dependence of the constant on $\|b\|_\infty$ is made explicit. Such dependence plays a crucial role in the distributional drift case discussed in Section \ref{Sec:pp}; yet, as far as we know, it has not been considered in the literature. Define
$$
\bP(t) := \mP \circ (X_t)^{-1}, \qquad \bP_n(t) := \mP \circ (X^n_t)^{-1},
$$
where $X^n$ is given by \eqref{eq:JN01-0}. The following theorem is our first main result.

\bt[Quantitative estimates: bounded drift]\label{thmEB}
Suppose that $T>0$, $\alpha \in (1,2)$, and $b \in L^\infty(\mR^d)$. Then
\begin{itemize}
\item[(i)] for any $\beta \in (0, (\alpha-1)/2)$ and $\delta \in (0, \alpha-1-\beta]$, there exists a constant $c$ depending only on $\Theta$, $\delta$, and $\|b\|_{\bB^{-\beta}_{\infty,\infty}}$ such that for all $n \in \mN$ and $t \in (0,T]$,
$$
\|\bP(t) - \bP_n(t)\|_{\rm var} \leq c \left(  \|b\|_\infty^{1+\delta} n^{-\delta} + \|b\|_\infty n^{-\delta/\alpha} + \|b\|_\infty^2 n^{-\frac{\alpha-1}{\alpha}} \right);
$$
\item[(ii)]  suppose that $\div b \in \bB^{-\beta}_{\infty,\infty}$, for any $\beta \in [(\alpha-1)/2, \alpha-1)$ and $\delta \in (0, \alpha-1-\beta]$, there exists a constant $c$ depending only on $\Theta$, $\delta$, $\|b\|_{\bB^{-\beta}_{\infty,\infty}}$, and $\|\div b\|_{\bB^{-\beta}_{\infty,\infty}}$ such that the same estimate as in (i) holds.
\end{itemize}
\et

\br
In particular, by setting $\delta = \alpha - 1 - \beta$, we obtain
\begin{align*}
\|\bP(t)-\bP_n(t) \|{\rm var} \leq c\left(  \|b\|_\infty^{\alpha-\beta} n^{-(\alpha-1-\beta)} + \|b\|_\infty n^{-\frac{\alpha-\beta-1}{\alpha}} + \|b\|_\infty^2 n^{-\frac{\alpha-1}{\alpha}} \right),
\end{align*} 
which matches the rate in \cite{SH24,FJM25} when $\beta = 0$, where the explicit dependence on $\|b\|_\infty$ in the constant was not provided in \cite{SH24,FJM25}.
\er

\subsection{Convergence rates for distributional drift}

Recall the mollified Euler's scheme \eqref{0608:00} for SDE \eqref{in:SDE} and denote
$$
\bP_{m,n}(t) := \mP \circ (X^{m,n}_t)^{-1}.
$$
Based on quantitative estimates for Euler's scheme with bounded drifts (see Theorem \ref{thmEB}) and the stability estimates (see Lemma \ref{thmSt}), we obtain our second main result: the weak convergence rates of the Euler--Maruyama scheme for SDEs driven by $\alpha$-stable processes with distributional drifts.

\bt[Weak convergence rates]\label{in:Main}
Assume that $T>0$, $\alpha \in (1,2)$,  and $\beta\in(0, \alpha-1 )$, and $m=n^\gamma$ with some $\gamma>0$.
\begin{itemize}
\item[(i)]  If $\beta\in (0,\tfrac{\alpha-1}{2})$ and $b\in \bB_{\infty,\infty}^{-\beta}$, then for any $\eps>0$, $\gamma\in(0,\frac{\alpha-1}{2\alpha\beta})$, and $\theta\in(\beta,\alpha-1-\beta)$, there is a constant $c>0$ depending only on $\Theta,\eps,\gamma,\theta,\phi,\|b\|_{\bB_{\infty,\infty}^{-\beta}}$ 
such that for any $n\in\mN$ and $t\in(0,T]$,
\begin{align*}
\|\bP(t)-\bP_{m,n}(t)\|_{\rm var} \leq c \left(  n^{  -\frac{\alpha-1}{\alpha}+\beta ( \gamma + \gamma \vee \frac{1}{\alpha} ) }+t^{\frac{\alpha-1-2\theta-\eps}{\alpha}}n^{-\gamma(\theta-\beta)} \right).
\end{align*}
\item[(ii)] If $\beta\in [\tfrac{\alpha-1}{2}, \alpha-1)$ and $b,\div b\in  \bB_{\infty,\infty}^{-\beta}$, then for any $\eps>0$ and $ \gamma\in(0, \frac{\alpha-1 -\beta}{\alpha\beta})$, there is a constant $c>0$ depending only on $\Theta,\eps, \gamma,\phi,\|b\|_{\bB_{\infty,\infty}^{-\beta}}, \|\div b\|_{\bB_{\infty,\infty}^{-\beta}}$ 
such that for any $n\in\mN$ and $t\in(0,T]$,
\begin{align*}
\|\bP(t)-\bP_{m,n}(t)\|_{\rm var} \leq c \left(     n^{  -\frac{\alpha-1}{\alpha}+\beta ( \gamma +  \frac{1}{\alpha} ) } 
+ n^{-\gamma(\alpha -1 -\beta) +\eps} \right ).
\end{align*}
\end{itemize}
\et

We illustrate our results by the following example.

\bexa 
If $\beta\in (0,\tfrac{\alpha-1}{2})$ and $b\in \bB_{\infty,\infty}^{-\beta}$, then
\begin{itemize}
\item[1)] for any small $\eps>0$, taking $\theta={(\alpha-1)}/{2}-\eps/\gamma$, we have that for any $n\in\mN$,
\begin{align*}
\sup_{t\in[0,T]}\|\bP_{m,n}(t)-\bP(t)\|_{\rm var} \lesssim  n^{-\frac{\alpha-1}{\alpha}  +\beta ( \gamma + \gamma \vee \frac{1}{\alpha} )  } +n^{-\gamma\frac{\alpha-1-2\beta}{2}+\eps},
\end{align*} 
which, when $\beta=0$ and $\gamma$ is taken sufficiently large, coincides with the rate $n^{-\frac{\alpha-1}{\alpha}}$ in \cite{SH24,FJM25} for the bounded-drift case;

\item[2)] for any small $\eps>0$, picking $\gamma=1/\alpha$ and $\theta=\alpha-1-\beta-\alpha\eps$, one sees that for any $n\in\mN$ and $t\in(0,T]$,
\begin{align*}
\|\bP_{m,n}(t)-\bP(t)\|_{\rm var} \lesssim   t^{-\frac{\alpha-1}{\alpha}} n^{-\frac{\alpha-1-2\beta}{\alpha}+\eps}.
\end{align*}
The rate of $n^{-\frac{\alpha-1-2\beta}{\alpha}}$ is natural considering the well-posedness condition $\beta \in (0, (\alpha-1)/2)$.
\end{itemize}
\eexa


\section{Preliminaries}\label{Pre}

\subsection{Besov spaces}\label{Sec:Be}

In this subsection, we introduce Besov spaces. Let $\sS(\mR^d)$ be the Schwartz space of all rapidly decreasing functions on $\mR^d$, and $\sS'(\mR^d)$ 
the dual space of $\sS(\mR^d)$ called Schwartz generalized function (or tempered distribution) space. Given $f\in\sS(\mR^d)$, 
the Fourier transform $\hat f$ and the inverse Fourier transform  $\check f$ are defined by
$$
\hat f(\xi) :=(2 \pi)^{-d/2}\int_{\mR^d} \e^{-i\xi\cdot x}f(x)\dif x, \quad\xi\in\mR^d,
$$
$$
\check f(x) :=(2 \pi)^{-d/2}\int_{\mR^d} \e^{i\xi\cdot x}f(\xi)\dif\xi, \quad x\in\mR^d.
$$
For every $f\in\sS'(\mR^d)$, the Fourier and the inverse  transforms are defined by
\begin{align*}
\<\hat{f},\varphi\>:=\<f,\hat{\varphi}\>,\qquad \<\check{f},\varphi\>:=\<f,\check{\varphi}\>, \ \  \forall\varphi\in\sS(\mR^d).
\end{align*}
Let $\chi:\mR^{d}\to[0,1]$ be a radial smooth function with
\begin{align*}
\chi(\xi)=
\begin{cases}
1, & \ \  |\xi|\leq 1,\\
0, &\ \ |\xi|>3/2.
\end{cases}
\end{align*}
For $\xi \in \mR^d$, define $\psi(\xi):=\chi(\xi)-\chi(2\xi)$ and for $j\in\mN_0$,
\begin{align*}
\psi_j(\xi){:=}\psi(2^{-j}\xi).
\end{align*}
Let $B_r := \{\xi\in \mR^d :  |\xi|\leq r\}$ for $r>0$. It is easy to see that $\psi\geq 0$,  supp$\psi\subset B_{3/2}/B_{1/2}$, and
\begin{align}\label{eq:SA00}
\chi(2\xi)+\sum_{j=0}^{k}\psi_j(\xi)=\chi(2^{-k}\xi)\to 1,\ \ \hbox{as}\ \ k\to\infty.
\end{align}
Since $\check \psi_j (y) = 2^{jd}\check \psi (2^j y), j \geq 0$, we have
\begin{align*}
\int_{\mR^d}|x|^\theta |\nabla^k\check \psi_j|(x) \dif x\leq {c} 2^{(k-\theta)j} , \ \ \theta>0,\ \ k\in\mN_0,
\end{align*}
where the constant $c$ is equal to $\int_{\mR^d}|x|^\theta |\nabla^k\check \psi|(x) \dif x$ and $\nabla^k$ stands for the $k$-order gradient. 
The block operators $\cR_j, j \geq 0$ are defined on $\sS'(\mR^d)$ by
\begin{align}\label{eq:Block}
\cR_j f (x):=(\psi_j\hat{f})^{\check\,}(x)=\check\psi_j* f (x) = 2^{jd} \int_{\mR^d} \check\psi(2^{j}y)f(x-y) \dif y,
\end{align}
and 
$
\cR_{-1}  f(x):= (\chi(2\cdot)\hat{f})^{\check\,}(x)=(\chi(2\cdot))\check{ }* f(x).
$
Then by \eqref{eq:SA00},
\begin{align}\label{eq:SA01}
f = \sum_{j \geq -1} \cR_j f.
\end{align}
 
Now we state the definitions of Besov spaces.
\bd[Besov spaces]\label{iBesov}
For every $s\in\mR$ and $p,q\in[1,\infty]$, the Besov space $\bB_{p,q}^s(\mR^d)$ is defined by
$$
\bB_{p,q}^s(\mR^d):=\Big\{f\in\sS'(\mR^d)\, \big| \, \|f\|_{\bB^s_{p,q}}:= \[ \sum_{j \geq -1}\left( 2^{s j} \|\cR_j f\|_{p} \right)^q \]^{1/q}  <\infty\Big\}.
$$
If $p=q=\infty$, it is in the sense
$$
\bB_{\infty,\infty}^s(\mR^d):=\Big\{f\in\sS'(\mR^d) \, \big| \, \|f\|_{\bB^s_{\infty,\infty}}:= \sup_{j \geq -1} 2^{s j} \|\cR_j f\|_{\infty} <\infty\Big\}.
$$
\ed 
 
Recall the following Bernstein's inequality (cf. \cite{BCD11}, Lemma 2.1).

\bl[Bernstein's inequality]\label{Bern}
For every $k\in\mN_0$, there is a constant $c=c(d,k)>0$ such that for all $j\ge-1$ and $1\leq p_1\leq p_2 \leq \infty$,
\begin{align*} 
\|\nabla^k\cR_j f\|_{p_2}  \leq c  2^{(k+ d (\frac{1}{p_1}-\frac{1}{p_2}))j}\|\cR_j f\|_{p_1}.
\end{align*}
In particular, for any $s \in \mR$ and $1\leq p, q \leq \infty$,
\begin{align}\label{S2:Bern}
\|\nabla^k f\|_{\bB^{s}_{p,q}} \leq c   \|f\|_{\bB^{s+k}_{p,q}}.
\end{align}
\el
\br
It is worth discussing here the equivalence between the Besov and H\"older  spaces, which will be used in various contexts in this paper without much explanation. For $s>0$, let $\bC^s(\mR^d)$ be the classical $s$-order H\"older space consisting of all measurable functions $f:\mR^d\to\mR$ with
\begin{align*}
\|f\|_{\bC^s}:=\sum_{j=0}^{[s]}\|\nabla^jf\|_\infty+[\nabla^{[s]}f]_{\bC^{s-[s]}}<\infty,
\end{align*}
where $[s]$ denotes  the largest integer less than or equal to $s$, and
\begin{align*}
\|f\|_\infty:=\sup_{x\in\mR^d}|f(x)|,\quad [f]_{\bC^\gamma}:=\sup_{h\in\mR^d}\frac{\|f(\cdot+h)-f(\cdot)\|_\infty}{|h|^\gamma},~\gamma\in(0,1).
\end{align*}
If $s>0$ and $s\notin\mN$, we have the  following equivalence between $\bB_{\infty,\infty}^s (\mR^d)$ and $ \bC^s  (\mR^d)$: (cf. \cite{Tr92})
\begin{align*}
\|f\|_{\bB_{\infty,\infty}^s}\asymp\|f\|_{\bC^s}.
\end{align*}
However, for any $n\in\mN_0$, we only have one side control that is
$
\|f\|_{\bB^n_{\infty,\infty}}\lesssim\|f\|_{\bC^n}.
$
\er

\br[Mollification in Besov spaces]
Let $\rho_m(x) :=m^{d}\rho (m  x)$, $m > 0$, be the mollifier for fixed $\rho \in C^{\infty}_c(\mR^d)$ being a smooth function with compact support and unit integral. Let $\beta\in\mR$ with $  \eps \in [0,1]$. It is easy to check that there is a constant $c>0$ such that for all $f\in\bB_{\infty,\infty}^{\beta+\eps}$ and $m \in\mN$,
\begin{align}\label{0725:new00}
    \|f-f_m\|_{\bB_{\infty,\infty}^{\beta}}\le c m^{-\eps}\|f\|_{\bB_{\infty,\infty}^{\beta+\eps}}.
\end{align}
\er

At the end of this subsection, we introduce the following   interpolation inequality (cf. \cite{BCD11}, Theorem 2.80).

\bl[Interpolation inequality]
Let $s_1,s_2\in \mR$ with $s_2>s_1$. For any $p\in [1,\infty]$ and $\theta \in (0,1)$, there is a constant $c=c(s_1,s_2,p)>0$ such that 
\begin{align*}
\|f\|_{\bB_{p,1}^{\theta s_1+ (1-\theta) s_2} } \leq c   \|f\|_{\bB^{s_1}_{p,\infty}}^\theta\|f\|_{\bB^{s_2}_{p,\infty}}^{1-\theta}.
\end{align*}
Furthermore, for any $ s_2> 0 > s_1$,
\begin{align}\label{InIn}
\|f\|_{\infty} \leq c   \|f\|_{\bB^{s_1}_{\infty,\infty}}^\theta\|f\|_{\bB^{s_2}_{\infty,\infty}}^{1-\theta},
\end{align}
where $\theta=s_2/(s_2-s_1)$.
\el

\subsection{$\alpha$-stable processes}

We call a $\sigma$-finite positive measure $\nu$ on $\mR^d$ a L\'evy measure if
\begin{align*}
\nu(\{0\})=0,\ \ \int_{\mR^d}\big(1\wedge|z|^2\big)\nu(\dif z)<+\infty.
\end{align*}
Fix $\alpha\in(0,2)$. Let $L^{(\alpha)}_t$ be a $d$-dimensional $\alpha$-stable process with L\'evy measure (or $\alpha$-stable measure) $\nu^{(\alpha)}$ defined as \eqref{CC02}. We say an $\alpha$-stable measure $\nu^{(\alpha)}$ is non-degenerate, if the assumption {\bf (ND)} holds. Note that for any $\gamma_2>\alpha>\gamma_1\geq 0$, 
\begin{align}\label{eq:BN01}
\int_{|z|\leq 1} |z|^{\gamma_2} \nu^{(\alpha)}(\dif z) + \int_{|z|>1} |z|^{\gamma_1}\nu^{(\alpha)}(\dif z)<\infty.
\end{align}
Let $N(\dif r, \dif z)$ be the associated Poisson random measure defined by
$$
N((0,t]\times A)  :=  \sum_{s\in(0,t]} \1_{A}(L_s^{(\alpha)} -  L^{(\alpha)}_{s-}),\ \ A\in \sB(\mR^d\setminus \{0\}) , t >0.
$$
By L\'evy-It\^o's decomposition (cf. \cite{Sa99}, Theorem 19.2), one sees that
\begin{align*}
L^{(\alpha)}_t =  \lim_{\eps \downarrow 0} \int_0^{t} \int_{\eps<|z|\leq 1} z \widetilde N(\dif r,\dif z) + \int_0^{t} \int_{|z|> 1} z   N(\dif r,\dif z),
\end{align*}
where $ \widetilde{N}(\dif r,\dif z):= N(\dif r,\dif z)-\nu^{(\alpha)}(\dif z)\dif r $ is the compensated Poisson random measure.

\subsection{Heat-kernel estimates}

Let $\alpha \in  (1, 2)$ and $L^{(\alpha)}$ be an $\alpha$-stable process having symmetric non-degenerate L\'evy measure $\nu^{(\alpha)}$. In this subsection, we start with the following time-inhomogeneous L\'evy process: for $0\leq t< \infty$,
\begin{align}\label{eq:L}
L^\sigma_t :=\int_0^t \sigma_r \dif L_r^{(\alpha)} = \int_0^{t} \int_{\mR^d} \sigma_r z \widetilde N(\dif r,\dif z),
\end{align}
where $\sigma : \mR_+ \to  \mM_{non}^{d}$ is  a bounded measurable function. Define
\begin{align}\label{eq:XM100}
P^\sigma_{s,t}f(x) := \mE f \left( x+\int_s^t \sigma_r \dif L_r^{(\alpha)} \right)
\end{align}
for all $f \in C_b^2(\mR^d)$. By It\^o's formula (cf. \cite{IW89}, Theorem 5.1 of Chapter II), one sees that
\begin{align}\label{eq:TE01}
\p_t P^\sigma_{s,t}f(x) = \sL^{(\alpha)}_{\sigma(t)}  P^\sigma_{s,t}f(x),
\end{align}
where 
\begin{align}\label{eq:a} 
\sL^{(\alpha)}_{\sigma(t)}  f(x):=\int_{\mR^d}\(f(x+\sigma(t) z)-f(x)- \sigma(t) z \cdot\nabla f(x)\)\nu^{(\alpha)}(\dif z).
\end{align} 

Below, we always make the following assumption in this subsection:

\medskip
\noindent
$\bf (\bH0)$ There is a constant $ \varkappa_0 >1$ such that
$$
 \varkappa_0^{-1}|\xi|\le|\sigma(t)\xi|\le  \varkappa_0|\xi|,\ \ \forall(t,\xi)\in\mR_+\times \mR^d.
$$

\medskip
\noindent
Under the assumptions $\bf (\bH0)$ and $\bf (ND)$, owing to  L\'evy-Khintchine's formula  (cf. \cite{Sa99}, Theorem 8.1) and \eqref{CC02}, for all $ |\xi|\geq1$, we have 
\begin{align*}
|\mE \e^{i\xi\cdot L^{\sigma}_{t} }|\leq&\exp\left( \int_0^t \int_{\mR^d}(\cos(\xi\cdot \sigma_s z)-1)\nu^{(\alpha)}(\dif z) \dif s \right)\\
\leq&\exp\left(- \int_0^t | \sigma_s^{\top}\xi|^\alpha \int_0^{\infty}\int_{\mS^{d-1}}\frac{1-\cos( \frac{\sigma_s^\top \xi}{|\sigma_s^\top  \xi |} \cdot r\theta)}{r^{1+\alpha}}\Sigma(\dif \theta)\dif r \dif s \right)
\leq \e^{-ct|\xi|^{\alpha}},
\end{align*}
where the constant $c>0$ depends only on $\alpha$, $\varkappa_0$, and $\Sigma(\mS^{d-1})$. Hence, by \cite{Sa99}, Proposition 28.1, the random variable $L_t^\sigma$ defined by \eqref{eq:L} admits a smooth density $p^\sigma(t,x)$ given by Fourier's inverse transform
$$
p^\sigma(t,x)=(2\pi)^{-d/2}\int_{\mR^d}\e^{-i x\cdot \xi}\mE \e^{i\xi\cdot L^\sigma_{t} }\dif \xi,\ \ \forall t>0,
$$
and the partial derivatives of $p^\sigma(t,\cdot)$ at any orders  tend to $0$ as $|x|\to \infty$.

The following integral-type estimate of  heat kernels is taken from \cite{CHZ20}, Lemma 3.2.

\bl
For each $0\leq s<t < \infty$, $p_{s,t}^\sigma(x)$  satisfies that for any $k\in\mN_0$ and $0\le\beta<\alpha$,  
\begin{align} \label{eq:HJ01}
\int_{\mR^d}|x|^\beta|\nabla^k p_{s,t}^\sigma (x)|\dif x \leq c  (t-s)^{-\frac{k-\beta}{\alpha}},
\end{align}
where  $c=c( \varkappa_0,k,d,\alpha,\beta)>0$.  
\el

From \eqref{eq:HJ01}, it is easy to check that for $f \in C_b^\infty(\mR^d)$,
\begin{align}\label{DH09}
\|\nabla^k P^\sigma_{s,t}f \|_\infty\lesssim (t-s)^{-\frac{k}{\alpha}}\|f\|_\infty, \quad  \text{for}~~k=0,1.
\end{align}
We also need the following heat kernel estimates in integral form with Littlewood-Paley’s decomposition, which is obtained in \cite{CHZ20}, Lemma 3.3 (see also \cite{HWW20}, Lemma 2.12). 

\bl
Suppose that $\bf (\bH0)$ holds with constant $\varkappa_0 >1$. Let $p_{s,t}^\sigma$ be the density of the random variable $L_t^\sigma - L_s^\sigma $. For any $n \in \mN_0$, and every $ \gamma\in[0,\alpha)$ and $ \vartheta\ge  \gamma$, there is a constant $c>0$ such that for all $0\leq s < t< \infty$ and $j\in \mN_0$,
\begin{align}\label{CC01}
\int_{\mR^d}|x|^{\gamma}|\nabla^n \cR_j p^\sigma_{s,t}(x)|\dif x \leq c  2^{(n-\vartheta)j} (t-s)^{-\frac{\vartheta}{\alpha}}\( (t-s)^{\frac{\gamma}{\alpha}} + 2^{-j\gamma}\),
\end{align}
where the block operators $\cR_j $ are defined by \eqref{eq:Block}. 
\el

We also need the following useful estimates.

\bl
Assume that $\alpha \in (1,2)$ and $T>0$. For $k=0,1$, there is a constant $c>0$ such that for all $0\leq u< s<t\le T$,
\begin{align}\label{eq:TE00}
\| \nabla^k \sL_{\sigma(t)}^{(\alpha)} P^\sigma_{s,t}f  \|_\infty \leq c (t-s)^{-\frac{k+\alpha}{\alpha}}\|f\|_\infty
\end{align}
and 
\begin{align}\label{DH10} 
\|\nabla^k P^\sigma_{u,t}f  -\nabla^k P^\sigma_{u,s}f \|_\infty\leq c \left[(s-u)^{-\frac{k}{\alpha}}\wedge ((s-u)^{-\frac{k+\alpha}{\alpha}}(t-s))\right]\|f\|_\infty.
\end{align}
\el

\begin{proof}
Observe that, under $\bf (\bH0)$, by \eqref{eq:a} and Bernstein's inequality,
\begin{align*}
\|\sL^{(\alpha)}_{\sigma(t)} \cR_j  h\|_\infty& \lesssim \int_{\mR^d} \( \[ |z| \| \nabla \cR_j  h \|_\infty\] \wedge \[|z|^2 \| \nabla^2 \cR_j  h \|_\infty\]\)\nu^{(\alpha)} (\dif z)\\
& \lesssim \| h \|_\infty \int_{\mR^d} \( |2^j z| \wedge |2^{j}z|^2\)\nu^{(\alpha)} (\dif z) \overset{\eqref{eq:BN01}}{\lesssim} 2^{\alpha j} \| h \|_\infty.
\end{align*}
Hence, by  \eqref{CC01} and Bernstein's inequality, we have
\begin{align*} 
\|\nabla^k \sL^{(\alpha)}_{\sigma(t)} P^\sigma_{s,t} f \|_\infty& \overset{\eqref{eq:SA01}}{ \lesssim}\sum_{j\ge-1}\|\sL^{(\alpha)}_{\sigma(t)} ( \nabla^k \cR_j   P^\sigma_{s,t} f)\|_\infty \overset{\eqref{eq:XM100}}{\lesssim} \sum_{j\ge-1}2^{(k+\alpha)j}\|\cR_j  p^\sigma_{s,t} \|_{1}\|f\|_\infty\\
&\lesssim \sum_{j\ge-1}2^{(k+\alpha)j}\left([2^{-(k+\alpha+1)j}(t-s)^{-\frac{k+\alpha+1}{\alpha}}]\wedge 1\right) \|f\|_\infty\\
& \lesssim (t-s)^{-\frac{k+\alpha}{\alpha}}\|f\|_\infty,
\end{align*}
where we used the following estimate in the last step: for any  $0<\beta<\gamma$ and $\lambda>0$,
\begin{align*}
\sum_{j\ge0} 2^{\beta j}\left([2^{-\gamma j}\lambda]\wedge 1\right)
& \le \lambda\wedge 1+ \int_0^\infty 2^{\beta s}\left([2^{-\gamma s}\lambda]\wedge 1\right)\dif s \\
&\lesssim \lambda\wedge 1 +\lambda^{\frac{\beta}{\gamma}}\int_{\lambda ^{-1/\gamma}}^\infty r^{  \beta-1}\left(r^{-\gamma}\wedge 1\right)\dif r \lesssim \lambda^{\frac{\beta}{\gamma}} .
\end{align*}
The first inequality \eqref{eq:TE00} follows. 

On the other hand, by \eqref{eq:TE01} and \eqref{eq:TE00}, for all $0\le s<t\le T$, we have
\begin{align*} 
|\nabla^k P^\sigma_{u,t} f (x)-\nabla^k P^\sigma_{u,s} f (x) |&  =\left|\int_s^t \nabla^k \p_r P^\sigma_{u,r} f(x) \dif r\right| =\left|\int_s^t \nabla^k\sL_{\sigma(r)}^{(\alpha)} P^\sigma_{u,r} (x) \dif r\right|\\
&\lesssim\|f\|_\infty\int_s^t
(r-u)^{-\frac{k+\alpha}{\alpha}}\dif r\\
& \lesssim  (s-u)^{-\frac{k+\alpha}{\alpha}}(t-s)\|f\|_\infty ,
\end{align*}
which, combining with \eqref{DH09}, deduces the desired result \eqref{DH10}.
\end{proof}

\section{Weak convergence rates}\label{sec:proof} 

\subsection{Euler's scheme for SDE with bounded drift}
Fix $T>0$.  In this subsection, we assume $b(x)$ belongs to $L^\infty(\mR^d)$ and consider the following SDE:
\begin{align}\label{in:SDE-b}
X_t=x+\int_0^tb(X_s)\dif s+ L_t^{(\alpha)},
\end{align}
and its Euler scheme: $X_0^n =X_0 = x$,
\begin{align}\label{eq:JN01}
X^n_t=x+\int_0^tb(X^n_{\pi_n(s)})\dif s+ L_t^{(\alpha)},
\end{align}
where  $n \in \mN$, and $\pi_n(t):=k/n$ for $t\in[k/n,(k+1)/n)$ with $k=0,1,2,...., \lfloor nT \rfloor$. Note that, for any $p \in (0,\alpha)$, by Lemma 2.10 in \cite{HW23},
\begin{align}\label{eq:mo}
\mE [|X^n_{r}-X^n_{\pi_n(r)}|^{p}]& \le \mE \(\|b\|_{\infty}n^{-1}+ |L^{(\alpha)}_r-L^{(\alpha)}_{\pi_n(r)}|\)^p\nonumber\\
& \leq (2^{p-1} \vee 1) \( \|b\|_\infty^pn^{-p}+ \mE  [|L^{(\alpha)}_r-L^{(\alpha)}_{\pi_n(r)}|^p] \)\nonumber\\ 
& \lesssim   \|b\|_\infty^pn^{-p}+ n^{-p/\alpha},
\end{align}
where the implicit constant in the inequality only depends on $d,\alpha,p,T$.

Now we are in a position to give 

\begin{proof}[Proof of Theorem \ref{thmEB}]

It suffices to estimate 
$$
\left | \mE\varphi(X^n_{t})-\mE\varphi(X_{t})\right |
$$
for any $\varphi \in C^\infty_b(\mR^d)$. The key ingredient of the proof is the Itô–Tanaka trick.

\medskip
\noindent
\medskip
\textbf{(Step 1)}  In this step, we prepare some estimates of PDEs for later use. Considering the following backward PDE with the terminal condition $\varphi\in C^\infty_b(\mR^d)$:
\begin{align}\label{eq:SZ00}
\p_su^t+\sL^{(\alpha)} u^t+b\cdot\nabla u^t=0,\quad u^t_t=\varphi,
\end{align}
where $u^t$ is the shifted function $u^t(s,x):=u(t-s,x)$ with $0\leq s<t\leq T$, and $\sL^{(\alpha)}$ is the infinitesimal  generator of $L_t^{(\alpha)}$ (see \eqref{eq:a}). It follows from Lemma 5.1 of \cite{HW23}  that for any $\beta\in (0,(\alpha-1)/2)$ (resp. $\beta \in [\frac{\alpha-1}{2}, \alpha-1)$) and $\delta\in[0,\alpha-\beta]$,
\begin{align}\label{DH12}
\|u^t(s)\|_{\bB_{\infty,\infty}^\delta}\lesssim (t-s)^{-\frac{\delta}{\alpha}}\|\varphi\|_\infty,
\end{align}
where the implicit constant in the above inequality only depends on $d,\alpha, T,\delta,\beta$, and $\|b\|_{ \bB^{-\beta}_{\infty,\infty}}$ (resp. $\|b\|_{ \bB^{-\beta}_{\infty,\infty}}, \|\div b\|_{ \bB^{-\beta}_{\infty,\infty}}$). Moreover, observe that, by interpolation inequality \eqref{InIn} and Bernstein's inequality \eqref{S2:Bern}, 
\begin{align}
\|\nabla u^t(s)\|_\infty & \lesssim \|\nabla u^t(s)\|^{1/2}_{\bB_{\infty,\infty}^{-(\alpha-\beta)+1}} \|\nabla u^t(s)\|^{1/2}_{\bB_{\infty,\infty}^{(\alpha-\beta)-1}} \nonumber\\
& \lesssim \| u^t(s)\|^{1/2}_{\bB_{\infty,\infty}^{2-(\alpha-\beta) }} \|  u^t(s)\|^{1/2}_{\bB_{\infty,\infty}^{\alpha-\beta }}\nonumber \\
&  \stackrel{\eqref{DH12}}{\lesssim} (t-s)^{-\frac{1}{\alpha}}. \label{eq:PL00}
\end{align}

\medskip
\noindent
\textbf{(Step 2)} In this step, we apply Itô's formula to rewrite $\mathbb{E}\varphi(X^n_t) - \mathbb{E}\varphi(X_t)$.  Adopting It\^o's formula (cf. \cite[Theorem 5.1 of Chapter II]{IW89}) to $u^t(s,X_s)$ and $u^t(s, X^n_s)$, one sees that,
\begin{align*}
u^t(s, X_s) & - u^t(0, x) =    \int_{0}^s (\partial_r u^t) (r, X_r)\dif r + \int_0^s b(X_{r})\cdot \nabla  u^t (r, X_r) \dif r\\
& + \int_0^s \int_{\mR^d} \left(u^t(r, X_{r-}+z)  - u^t(r, X_{r-}) \right) \widetilde N(\dif r, \dif z)\\
& + \int_0^s \int_{\mR^d} \Big(u^t (r, X_r+z) - u^t (r, X_r)   - z \cdot \nabla u^t (r, X_r) \Big) \nu^{(\alpha)}(\dif z)\dif r,
\end{align*}
and
\begin{align*}
u^t(s, X^n_s)&  - u^t(0, x) = \int_{0}^s (\partial_r u^t) (r, X^n_r)\dif r + \int_0^s b(X^n_{\pi_n(r)})\cdot \nabla  u^t (r, X^n_r) \dif r\\
& + \int_0^s \int_{\mR^d} \left(u^t(r, X^n_{r-}+z)  - u^t(r, X^n_{r-}) \right) \widetilde N(\dif r, \dif z)\\
& + \int_0^s \int_{\mR^d} \Big(u^t (r, X_r^n+z) - u^t (r, X_r^n)   - z \cdot \nabla u^t (r, X_r^n) \Big) \nu^{(\alpha)}(\dif z)\dif r,
\end{align*}
where the third terms on the right-hand side of the above equalities are martingales. Then by \eqref{eq:SZ00}, it is easy to check that
\begin{align*}
u^t(s, X_s)  - u^t(0, x) =     \int_0^s \int_{\mR^d} \left(u^t(r, X_{r-}+z)  - u^t(r, X_{r-}) \right) \widetilde N(\dif r, \dif z) ,
\end{align*}
and
\begin{align*}
u^t(s, X^n_s)  - u^t(0, x)   =  & \int_0^s \( b(X^n_{\pi_n(r)}) - b(X^n_{r}) \) \cdot \nabla  u^t (r, X^n_r) \dif r\\
& + \int_0^s \int_{\mR^d} \left(u^t(r, X^n_{r-}+z)  - u^t(r, X^n_{r-}) \right) \widetilde N(\dif r, \dif z) .\end{align*}
Furthermore,   taking $s=t$, we have
\begin{align*}\mE\varphi(X_{t})= \mE u^{t}(t,X_{t}) = u^t(0,x),
\end{align*}
and
\begin{align*}\mE\varphi(X^n_{t})= \mE u^{t}(t,X^n_{t}) = u^t(0,x)+\mE\int_0^t
\left(b(X^n_{\pi_n(r)})-b(X^n_{r})\right)\cdot\nabla u^t(r,X^n_{r})\dif r.
\end{align*}
Thus, we get
 \begin{align}\label{eq:YH00}
 \mE\varphi(X^n_{t})-\mE\varphi(X_{t}) 
 = \mE\int_0^t
\left(b(X^n_{\pi_n(r)})-b(X^n_{r})\right)\cdot\nabla u^t(r,X^n_{r})\dif r.
\end{align}

\medskip
\noindent
\medskip
\textbf{(Step 3)} Thanks to  \eqref{eq:YH00}, we have
 \begin{align*}
\mE\varphi(X^n_{t})-\mE\varphi(X_{t}) 
 = & \,\mE\int_0^t b (X^n_{\pi_n(r)} )\cdot\left(\nabla u^t (r,X^n_{r})-\nabla u^t(r,X^n_{\pi_n(r)})\right)\dif r\\
& +  \int_0^t\left[\mE \(b\cdot\nabla u^t(r)\)(X^n_{\pi_n(r)})-\mE \(b\cdot\nabla u^t(r)\)(X^n_{r}))\right]\dif r\\
=:&\, \sI_1(t)+\sI_2(t).
\end{align*}
Next, we estimate these two terms in turn.

\medskip
\noindent \textbf{(Step 3.1)}
For $\sI_1(t)$, by Bernstein's inequality \eqref{S2:Bern}, and \eqref{eq:mo}, one sees that for any $\delta \in (0, \alpha-1-\beta]$,
\begin{align*}
|\sI_1(t)| & \lesssim  \|b\|_\infty\int_0^t \|\nabla u^t(r)\|_{\bB_{\infty,\infty}^\delta} \mE|X^n_{\pi_n(r)}-X^n_{r}|^{\delta}\dif r\\
& \lesssim  \|b\|_\infty \int_0^t \|\  u^t(r)\|_{\bB_{\infty,\infty}^{1+\delta}} \( \|b\|_\infty^\delta n^{-\delta}+ n^{-\delta/\alpha} \)  \dif r\\
& \overset{\eqref{DH12}}{\lesssim}  \|b\|_\infty \( \|b\|_\infty^\delta n^{-\delta}+ n^{-\delta/\alpha} \)  \int_0^t(t-r)^{-\frac{1+\delta}{\alpha}}  \dif r\\
& =  \|b\|_\infty ( \|b\|_\infty^\delta n^{-\delta}+ n^{-\delta/\alpha}) t^{\frac{\alpha-1-\delta}{\alpha}}\int_0^1r^{\frac{\alpha-1-\delta}{\alpha} -1} \dif r.
\end{align*}
Consequently, we get that for each $\delta \in (0, \alpha-1-\beta]$,
\begin{align}\label{eq:RR01}
|\sI_1(t)| \lesssim  \|b\|_\infty^{1+\delta} n^{-\delta}+  \|b\|_\infty n^{-\delta/\alpha}, ~~ \text{for}~~ t \in [0,T].
\end{align}

\medskip
\noindent \textbf{(Step 3.2)}
As for $\sI_2(t)$,  the estimate of 
\begin{align}\label{eq:HJ00}
\left | \mE \left[ \( b\cdot\nabla u^t(r) \)(X^n_{\pi_n(r)}) \right] -\mE \left[ \( b\cdot\nabla u^t(r) \)(X^n_{r}) \right] \right |,
\end{align}
is the key ingredient. 

\medskip
\noindent (i)
Using the It\^o-Tanaka trick again, we consider the following equation:
\begin{align}\label{eq:AJ01}
\p_s w^r+\sL^{(\alpha)} w^r=0,\quad w^r(r)=f,
\end{align}
where $w^r(s,x):=w(r-s,x)$ is the shifted function with $0\leq s<r\leq T$ and $f\in C^\infty_b(\mR^d)$. Below, we will take $f = b\cdot\nabla u^t(r)$ in the step (ii).  Applying It\^o's formula (cf. Theorem 5.1 of Chapter II in \cite{IW89}) to $w^r(s,X^n_s)$ and by \eqref{eq:AJ01}, we have
\begin{align*}
w^r(t', X^n_{t'}) & - w^r(0, x) =    \int_0^{t'} b(X^n_{\pi_n(s)})\cdot\nabla w^r(s,X^n_s)\dif s\\
& +  \int_0^{t'} \int_{\mR^d} \left(w^r (s, X^n_{s-}+z)  - w^r (s, X^n_{s-}) \right) \widetilde N(\dif s, \dif z),
\end{align*}
which implies that
\begin{align*}
\mE f(X^n_r) = \mE w^r(r,X^n_r)=w^r(0,x)+\mE\int_0^r b(X^n_{\pi_n(s)})\cdot\nabla w^r(s,X^n_s)\dif s.
\end{align*}
Hence, for any $r_2>r_1$,
\begin{align*}
\begin{split}
\left|\mE f(X^n_{r_2})-\mE f(X^n_{r_1})\right| \leq \, & |w(r_2,x) - w(r_1,x)|  + \left|\mE\int_{r_1}^{r_2} b(X^n_{\pi_n(s)})\cdot\nabla w(r_2-s,X^n_s)\dif s\right|\\
&+ \left|\mE\int_0^{r_1} b(X^n_{\pi_n(s)})\cdot \(\nabla w(r_2-s,X^n_s) - \nabla w(r_1-s,X^n_s) \)\dif s\right|\\
\leq\, & \| w(r_2)-w(r_1)\|_\infty+\|b\|_\infty\int_{r_1}^{r_2}\|\nabla w(r_2-s)\|_\infty\dif s\\
&+\|b\|_\infty\int_0^{r_1}\|\nabla w(r_2-s)-\nabla w(r_1-s)\|_\infty\dif s\\
 := \, &  \sA_{2,1} + \sA_{2,2} + \sA_{2,3}.
\end{split}
\end{align*}
By \eqref{DH10}, one sees that
\begin{align}\label{eq:RR02}
\sA_{2,1} \lesssim \|f\|_\infty
 \left[1\wedge ({r_1}^{-1}(r_2-{r_1}))\right].
 \end{align}
Based on  \eqref{DH09}, we obtain that 
\begin{align}\label{eq:RR03}
 \sA_{2,2} \lesssim
 \|b\|_\infty \int_{r_1}^{r_2}  (r_2-s )^{-\frac{1}{\alpha}} \dif s
 \lesssim 
  \|b\|_\infty\|f\|_\infty (r_2-{r_1})^{-\frac{1}{\alpha}+1}.
\end{align}
Using \eqref{DH10} again, we have that for all $0< r_1<r_2\le T$,
\begin{align}\label{eq:RR04}
 \sA_{2,3}& \lesssim \|b\|_\infty\|f\|_\infty \int_0^{r_1}\left[({r_1}-s)^{-\frac{1}{\alpha}}\wedge ((r_2-{r_1})({r_1}-s)^{-\frac{1+\alpha}{\alpha}})\right]\dif s \nonumber\\
 & =  \|b\|_\infty\|f\|_\infty \int_0^{r_1} s^{-\frac{1}{\alpha}} \left[1 \wedge \((r_2-{r_1})s^{-1}\)\right]\dif s \nonumber\\
 &\lesssim \|b\|_\infty \|f\|_\infty  \left[ \int_0^{r_2-r_1} s^{-\frac{1}{\alpha}} \dif s+ (r_2-{r_1}) \int_{(r_2-r_1)\wedge r_1}^{r_1}   s^{-\frac{1+\alpha}{\alpha}} \dif s \right ] \nonumber\\
 & \lesssim  \|b\|_\infty \|f\|_\infty (r_2-{r_1})^{-\frac{1}{\alpha}+1} .
\end{align}
Combining the estimates \eqref{eq:RR02}-\eqref{eq:RR04}, we obtain that  for all $0< r_1<r_2\le T$,
\begin{align}\label{DH11}
\left|\mE f(X^n_{r_2})-\mE f(X^n_{r_1})\right|
\lesssim \|f\|_\infty\left(\left[1\wedge ({r_1}^{-1}(r_2-{r_1}))\right]+\|b\|_\infty(r_2-{r_1})^{-\frac{1}{\alpha}+1}\right).
\end{align} 

 \medskip
\noindent (ii)
Next, we substitute $(f(\cdot), r_1, r_2)$ in \eqref{DH11} with 
$$
((b \cdot \nabla u^t(r))(\cdot), \pi_n(r), r)
$$
to estimate \eqref{eq:HJ00} and then $\sI_2(t)$. Observe that this substitution is only justified when $r \ge 1/n$, since \eqref{DH11} requires $0 < r_1 < r_2 \le T$. 

 \medskip
 
We therefore distinguish two cases. 

 \medskip
\noindent $\bullet$
If $t< 1/n$, then 
\begin{align*}
|\sI_2(t)|
&  \lesssim \|b\|_\infty \int_0^t \|\nabla u^t(r)\|_\infty \dif r 
 \stackrel{\eqref{eq:PL00}}{\lesssim} \|b\|_\infty \int_0^t (t-r)^{-\frac1\alpha}\dif r \\
& \lesssim \|b\|_\infty t^{\frac{\alpha-1}{\alpha}}
\lesssim \|b\|_\infty n^{-\frac{\alpha-1}{\alpha}}.
\end{align*}

 \medskip
\noindent $\bullet$
If $t\ge 1/n$, then we split
\begin{align*}
|\sI_2(t)|
 \leq \, & \left(  \int_0^{1/n} + \int_{1/n}^{t} \right) \left|\mE  \big[ ( b\cdot\nabla u^t(r) )(X^n_{\pi_n(r)})  \big] -\mE \left[ ( b\cdot\nabla u^t(r) )(X^n_{r}) \right] \right| \dif r\\
 =: \, &   \sI_{2,0}(t)+  \sI_{2,1}(t).
\end{align*}
For $\sI_{2,0}(t)$ on $t\in [ 1/n,T]$, using the trivial bound, we get
\begin{align*}
\sI_{2,0}(t)
&\lesssim \int_0^{1/n}\|b\cdot\nabla u^t(r)\|_\infty \dif r \lesssim \|b\|_\infty\int_0^{1/n}\|\nabla u^t(r)\|_\infty\dif r\\
&\stackrel{\eqref{eq:PL00}}{\lesssim} \|b\|_\infty\int_0^{1/n}(t-r)^{-\frac1\alpha}\dif r \le \|b\|_\infty\int_{0}^{1/n}\left(\frac1n - r\right)^{-\frac1\alpha}\dif r\\
&\lesssim \|b\|_\infty\int_0^{1/n}s^{-\frac1\alpha}\dif s
\lesssim \|b\|_\infty n^{-\frac{\alpha-1}{\alpha}}.
\end{align*}
For $ \sI_{2,1}(t)$ on $t \in  [ 1/n,T]$, since the integration variable now satisfies $r \ge 1/n$, we have $\pi_n(r)>0$, and thus \eqref{DH11} applies with $(f (\cdot), r_1,r_2)=((b \cdot \nabla u^t(r))(\cdot), \pi_n(r),r)$. In this case, we infer that since $\frac{r}{\pi_n(r)} \leq 2$ (for $r>1/n$) and $|\pi_n(r) -r|\leq 1/n$,
\begin{align*}
 \sI_{2,1}(t)
& \lesssim    \int_{1/n}^t \|b \cdot \nabla u^t(r)\|_\infty \Bigg\{ \left[1\wedge \((\pi_n(r))^{-1}(r- \pi_n(r))\)\right]+\|b\|_{\infty}(r- \pi_n(r))^{\frac{\alpha-1}{\alpha}}\Bigg\}\dif r\\
& \lesssim  \|b\|_\infty \int_{1/n}^t \|\nabla u^t(r)\|_\infty \left(\left[1\wedge ((nr)^{-1})\right]+\|b\|_{\infty}n^{-\frac{\alpha-1}{\alpha}}\right)\dif r\\
& \overset{\eqref{eq:PL00}}{\lesssim} n^{-\frac{\alpha-1}{\alpha}}\|b\|_\infty\left( \int_0^t(t-r)^{-\frac{1}{\alpha}}r^{-\frac{\alpha-1}{\alpha}}\dif r+\|b\|_\infty   \right)\\
& \lesssim  n^{-\frac{\alpha-1}{\alpha}}(\|b\|_\infty+\|b\|_\infty^2),
\end{align*}
where we used the fact $1\wedge a^{-1} \leq a^{-x}$ for $\forall a>0,x \in [0,1]$ in the third inequality, and the definition of the Beta function \eqref{eq:Beta} in the last inequality. 

 \medskip
 
Consequently, we have
\begin{align*}
|\sI_2(t)|  \lesssim  n^{-\frac{\alpha-1}{\alpha}}(\|b\|_\infty+\|b\|_\infty^2), ~~ 
\text{for}~~t \in [0,T],
\end{align*}
which together with \eqref{eq:RR01} (the estimate of $|\sI_1(t)|$ on $[0,T]$) derives that for $t \in [0,T]$ and $\delta \in (0, \alpha-1-\beta]$, 
\begin{align*}
& \quad | \mE \varphi(X^n_{t})-\mE \varphi(X_{t}) |  \leq  |\sI_1(t)| + |\sI_2(t)|\\
&  \lesssim \|b\|_\infty^{1+\delta} n^{-\delta}+  \|b\|_\infty n^{-\delta/\alpha} + n^{-\frac{\alpha-1}{\alpha}}(\|b\|_\infty+\|b\|_\infty^2)\\
& \leq \|b\|_\infty^{1+\delta} n^{-\delta}+  \|b\|_\infty n^{-\delta/\alpha} +\|b\|_\infty^2 n^{-\frac{\alpha-1}{\alpha}} .
\end{align*}
This yields the desired estimates.
\end{proof}

\subsection{Euler's scheme for SDE with distributional drift}

To prove Theorem \ref{in:Main}, we need the following stability estimates taken from  \cite{HW23}. 

\bl[Stability estimates]\label{thmSt}
 Let $T>0$, $\alpha\in(1,2)$ and $\beta\in [0,\alpha-1)$.
Assume that  $X^1$ and $X^2$ are weak solutions to SDE \eqref{in:SDE} with drift $b=b_1 \in  \bB_{\infty,\infty}^{-\beta}$ and $b=b_2\in  \bB_{\infty,\infty}^{-\beta}$, respectively. Denote the time marginal law of $X^i$ by $\bP_i(t)$, $i=1,2$. Then,
\begin{itemize}
\item[(i)] when $\beta<\frac{\alpha-1}{2}$, for any  $\theta\in[\beta,\alpha-1-\beta)$ and $\eps>0$, there is a constant $c=c (\Theta,\theta,\eps,\|b_1\|_{ \bB_{\infty,\infty}^{-\beta}})>0$ such that for any $t\in (0,T]$,
\begin{align*}
\|\bP_1(t)-\bP_2(t)\|_{\rm var}  \leq c   t^{\frac{\alpha-1-2\theta-\eps}{\alpha} }\|b_1-b_2\|_{ \bB_{\infty,\infty}^{-\theta}};
\end{align*}
\item[(ii)] when $\beta\ge\frac{\alpha-1}{2}$, for any  $\theta\in[\beta,\alpha-1)$ and $\eps>0$, there is a constant $c>0$ depending on $\Theta,\theta,\eps,\|b_1\|_{\bB_{\infty,\infty}^{-\beta}},\|\div b_1\|_{\bB_{\infty,\infty}^{-\beta}}$, such that for any $t\in (0,T]$,
\begin{align*}
\|\bP_1(t)-\bP_2(t)\|_{\rm var}  \leq c   t^{\frac{\alpha-1-\theta-\eps}{\alpha} }\left(\|b_1-b_2\|_{ \bB_{\infty,\infty}^{-\theta}}+\|\div b_1-\div b_2\|_{\bB_{\infty,\infty}^{-\theta}}\right).
\end{align*}

\end{itemize}

\el

\br 
The stability result clearly indicates that the weak solution obtained in Proposition \ref{thm:G-mart} is independent of the specific choice of mollifier functions $\phi_m$.
\er

Now we are in a position to give 

\begin{proof}[Proof of  Theorem \ref{in:Main}]
Let $X^m_t$ be the solution to the following classical SDE:
\begin{align*}
X_t^m=x+\int_0^tb_m(X^m_s)\dif s+L_t^{(\alpha)},
\end{align*}
where $b_m$ is defined by \eqref{BN1}. Denote $\bP_m(t):=\mP\circ(X_t^m)^{-1}$. Note that 
$$
\|b_m\|_\infty\lesssim m^\beta\|b\|_{\bB_{\infty,\infty}^{-\beta}}=n^{\beta\gamma} \|b\|_{\bB_{\infty,\infty}^{-\beta}}.
$$

\noindent
\textbf{(Step 1)}
Applying Theorem \ref{thmEB} with $\delta = \alpha -1 - \beta$, for any $ \gamma>0$, we get that
\begin{align*}
\|\bP_{m,n}(t)-\bP_m(t)\|_{\rm{var}} 
& \lesssim  \|b_m\|_\infty^{\alpha-\beta}n^{-(\alpha-1-\beta)}+\|b_m\|_\infty n^{-\frac{\alpha-1-\beta}{\alpha}}+ \|b_m\|_\infty^2 n^{-\frac{\alpha-1}{\alpha}}\\
&\lesssim  n^{-(\alpha-1-\beta)+\beta\gamma(\alpha-\beta)}+n^{-\frac{\alpha-\beta-1}{\alpha}+\beta\gamma}+n^{-\frac{\alpha-1}{\alpha}+2\beta\gamma}\\
& \lesssim   n^{-\frac{\alpha-1}{\alpha}+   \beta ( \gamma + \gamma \vee \frac{1}{\alpha} )}.
\end{align*} 
                                                                                                                                                                                                                                                                                                                                                                                                                                                                                                                                                                                                                                                                                                       Notice that 
\begin{align}\label{eq:YR00}
-\frac{\alpha-1}{\alpha} + \beta(\gamma + \gamma \vee \tfrac{1}{\alpha}) < 0
\end{align}
if and only if
$$
\gamma \in \Big(~0,~ \frac{1}{\alpha}\big[\big(\tfrac{\alpha-1}{\beta} - 1\big) \wedge 1 \big] ~ \Big) \cup \Big[~ \dfrac{1}{\alpha},~ \dfrac{\alpha-1}{2\alpha \beta } ~\Big).
$$
Thus, if $\beta< \frac{\alpha-1}{2}$, then
$$
\frac{\alpha-1}{\beta} - 1 > \frac{\alpha-1}{2\beta} \geq 1 \Rightarrow 0< \gamma <  \frac{\alpha-1}{2\alpha\beta}  \Rightarrow \eqref{eq:YR00};
$$
if $\beta \in [ \frac{\alpha-1}{2}, \alpha-1)$, then   
$$
 \frac{\alpha-1}{\beta} - 1 \leq \frac{\alpha-1}{2\beta}\leq 1  \Rightarrow 0<\gamma < \frac{\alpha-1 -\beta}{\alpha\beta}\leq  1/\alpha  \Rightarrow \eqref{eq:YR00}.
$$
Consequently, we have
\begin{align*}
\|\bP_{m,n}(t)-\bP_m(t)\|_{\rm{var}} 
& \lesssim \begin{cases}
n^{-\frac{\alpha-1}{\alpha}+   \beta ( \gamma + \gamma \vee \frac{1}{\alpha} )} ,&   \text{if}~~\beta \in (0, \frac{\alpha-1}{2}), \gamma \in (0, \frac{\alpha-1}{2\alpha\beta}), \\
n^{-\frac{\alpha-1}{\alpha}+   \beta ( \gamma +  \frac{1}{\alpha} )}, &   \text{if}~~\beta \in [ \frac{\alpha-1}{2}, \alpha-1), \gamma\in (0, \frac{\alpha-1 -\beta}{\alpha\beta}).
\end{cases}
\end{align*}

\medskip
\noindent 
\textbf{(Step 2)}
Moreover, according to the stability estimates Lemma \ref{thmSt} and taking $m = n^{\gamma}$, one sees that

\medskip\noindent
{\textbf{(i)}} when $\beta < \frac{\alpha-1}{2}$, for any $\eps>0$ and $\theta\in(\beta,\alpha-1-\beta)$, 
\begin{align*}
\|\bP_m(t)-\bP(t)\|_{\rm{var}}&
\lesssim t^{\frac{\alpha-1-2\theta-\eps}{\alpha}}\|b-b_m\|_{\bB_{\infty,\infty}^{-\theta}}\nonumber\\
& \lesssim t^{\frac{\alpha-1-2\theta-\eps}{\alpha}}m^{-{(\theta-\beta)}} \|b\|_{\bB_{\infty,\infty}^{-\beta}}\\
&  \lesssim  t^{\frac{\alpha-1-2\theta-\eps}{\alpha}}n^{-\gamma(\theta-\beta)} ,
\end{align*}
where we used \eqref{0725:new00} in the second inequality;

\medskip\noindent
{\textbf{(ii)}}  when $\beta \geq \frac{\alpha-1}{2}$, for any $\eps>0$ and $\theta\in(\beta,\alpha-1 )$,  
\begin{align*}
\|\bP_m(t)-\bP(t)\|_{\rm{var}}&
\lesssim t^{\frac{\alpha-1- \theta-\eps}{\alpha}} \left( \|b-b_m\|_{\bB_{\infty,\infty}^{-\theta}}+  \|\div b-\div b_m\|_{\bB_{\infty,\infty}^{-\theta}} \right)\nonumber\\
&  \lesssim t^{\frac{\alpha-1- \theta-\eps}{\alpha}}  \left( m^{-(\theta-\beta)} \| b \|_{\bB_{\infty,\infty}^{-\beta}}+  m^{-(\theta-\beta)} \|\div  b \|_{\bB_{\infty,\infty}^{-\beta}}\right)  \\ 
& \lesssim  t^{\frac{\alpha-1- \theta-\eps}{\alpha}} n^{-\gamma(\theta-\beta)} .
\end{align*} 
Furthermore, for $0<\gamma < \frac{\alpha-1 -\beta}{\alpha\beta} \leq 1$, taking $\theta = \alpha -1 - \eps/\gamma$ with small enough $\eps>0$, we have that
$$
\|\bP_m(t)-\bP(t)\|_{\rm{var}} \lesssim n ^{-\gamma(\alpha - 1- \beta) + \eps}.
$$

Finally, combining the above calculations and observing 
\begin{align*}
\|\bP_{m,n}(t)-\bP(t)\|_{\rm{var}} 
&\leq \|\bP_{m,n}(t)-\bP_m(t)\|_{\rm{var}} + \|\bP_m(t)-\bP(t)\|_{\rm{var}},
\end{align*}
we get the desired result.
\end{proof}

\subsection*{Acknowledgements}


 Mingyan Wu is supported by the National Natural Science Foundation of China (Grant No. 12201227).

\begin{appendix}

\end{appendix}

\addcontentsline{toc}{section}{References}


\begin{center}
\scshape References
\end{center}



\begin{biblist}


\bib{ABM18}{article}{
   author={Athreya, Siva},
   author={Butkovsky, Oleg},
   author={Mytnik, Leonid},
   title={Strong existence and uniqueness for stable stochastic differential
   equations with distributional drift},
   journal={Ann. Probab.},
   volume={48},
   date={2020},
   number={1},
   pages={178--210},
   issn={0091-1798},
   review={\MR{4079434}},
   doi={10.1214/19-AOP1358},
}
\bib{BCD11}{book}{
      author={Bahouri, Hajer},
      author={Chemin, Jean-Yves},
      author={Danchin, Rapha\"{e}l},
       title={Fourier analysis and nonlinear partial differential equations},
      series={Grundlehren der Mathematischen Wissenschaften [Fundamental
  Principles of Mathematical Sciences]},
   publisher={Springer, Heidelberg},
        date={2011},
      volume={343},
        ISBN={978-3-642-16829-1},
         url={https://doi.org/10.1007/978-3-642-16830-7},
      review={\MR{2768550}},
}  
   
\bib{CC18}{article}{
   author={Cannizzaro, G.},
   author={Chouk, K.},
   title={Multidimensional SDEs with singular drift and universal
   construction of the polymer measure with white noise potential},
   journal={Ann. Probab.},
   volume={46},
   date={2018},
   number={3},
   pages={1710--1763},
   issn={0091-1798},
   review={\MR{3785598}},
   doi={10.1214/17-AOP1213},
}
\bib{CHT22}{article}{
   author={Cannizzaro, Giuseppe},
   author={Haunschmid-Sibitz, Levi},
   author={Toninelli, Fabio},
   title={$\sqrt{\log t}$-superdiffusivity for a Brownian particle in the
   curl of the 2D GFF},
   journal={Ann. Probab.},
   volume={50},
   date={2022},
   number={6},
   pages={2475--2498},
   issn={0091-1798},
   review={\MR{4499841}},
   doi={10.1214/22-aop1589},
}

 
\bib{Ca22}{article}{
   author={Cavallazzi, Thomas},
   title={Quantitative weak propagation of chaos for McKean-Vlasov SDEs
   driven by stable processes},
   language={English, with English and French summaries},
   journal={Ann. Inst. Henri Poincar\'e{} Probab. Stat.},
   volume={61},
   date={2025},
   number={3},
   pages={1662--1764},
   issn={0246-0203},
   review={\MR{4947166}},
   doi={10.1214/24-aihp1475},
    eprint={2212.01079},
}
 
 \bib{CIP25}{article}{
   author={Chaparro J\'aquez, Luis Mario},
   author={Issoglio, Elena},
   author={Palczewski, Jan},
   title={convergence rates of numerical scheme for SDEs with a
   distributional drift in Besov space},
   journal={ESAIM Math. Model. Numer. Anal.},
   volume={59},
   date={2025},
   number={5},
   pages={2717--2738},
   issn={2822-7840},
   review={\MR{4969855}},
   doi={10.1051/m2an/2025064},
}
\bib{CMOW23}{article}{
author={Chatzigeorgiou, Georgiana},
author={Morfe, Peter},
author={Otto, Felix},
author={Wang, Lihan},
title={The Gaussian free-field as a stream function: asymptotics of effective diffusivity in infra-red cut-off},
eprint={2212.14244}
}

\bib{CM19}{article}{
   author={Chaudru de Raynal, Paul-\'Eric},
   author={Menozzi, St\'ephane},
   title={On multidimensional stable-driven stochastic differential
   equations with Besov drift},
   journal={Electron. J. Probab.},
   volume={27},
   date={2022},
   pages={Paper No. 163, 52},
   review={\MR{4525442}},
   doi={10.1214/22-ejp864},
}

\bib{CHZ20}{article}{
      author={Chen, Zhen-Qing},
      author={Hao, Zimo},
      author={Zhang, Xicheng},
       title={H\"{o}lder regularity and gradient estimates for {SDE}s driven by
  cylindrical {$\alpha$}-stable processes},
        date={2020},
     journal={Electron. J. Probab.},
      volume={25},
       pages={Paper No. 137, 23},
         url={https://doi.org/10.1214/20-ejp542},
   review={\MR{4179301}},
}

 
\bib{DGI22}{article}{author={De Angelis, Tiziano},
  author={Germain, Maximilien},
   author={Issoglio, Elena},
   title={A numerical scheme for stochastic differential equations with
   distributional drift},
journal={Stochastic Process. Appl.},
   volume={154},
   date={2022},
   pages={55--90},
  issn={0304-4149},
review={\MR{4490482}},
   doi={10.1016/j.spa.2022.09.003},
}
\bib{DD16}{article}{
   author={Delarue, F.},
   author={Diel, R.},
   title={Rough paths and 1d SDE with a time dependent distributional drift:
   application to polymers},
   journal={Probab. Theory Related Fields},
   volume={165},
   date={2016},
   number={1-2},
   pages={1--63},
   issn={0178-8051},
   review={\MR{3500267}},
   doi={10.1007/s00440-015-0626-8},
}

\bib{FJM25}{article}{     
	author={Fitoussi, Mathis},     
	author={Jourdain, Benjamin},     
	author={Menozzi, Stéphane},
       title={Weak well-posedness and weak discretization error for stable-driven SDEs with Lebesgue drift},   
       journal={IMA Journal of Numerical Analysis},
       year={2025},
       eprint={2405.08378},
       doi={10.1093/imanum/draf079}
}

 \bib{GHR25}{article}{
 author={Gouden{\`e}ge, Ludovic},
 author={Haress, El Mehdi},
 author={Richard, Alexandre},
 issn={0304-4149},
 doi={10.1016/j.spa.2024.104533},
 review={Zbl 1557.60160},
 title={Numerical approximation of SDEs with fractional noise and distributional drift},
 journal={Stochastic Processes and their Applications},
 volume={181},
 pages={38},
 note={Id/No 104533},
 date={2025},
 publisher={Elsevier (North-Holland), Amsterdam},
}


\bib{GK}{article}{
   author={Gy\"ongy, Istv\'an},
   author={Krylov, Nicolai},
   title={Existence of strong solutions for It\^o's stochastic equations via
   approximations},
   journal={Probab. Theory Related Fields},
   volume={105},
   date={1996},
   number={2},
   pages={143--158},
   issn={0178-8051},
   review={\MR{1392450}},
   doi={10.1007/BF01203833},
}


\bib{HW23}{article}{
      author={Hao, Zimo},
   author={Wu, Mingyan},
       title={SDE driven by multiplicative cylindrical $\alpha$-stable noise with distributional drift},
    eprint={2305.18139},
}
 
\bib{HWW20}{article}{
   author={Hao, Zimo},
   author={Wang, Zhen},
   author={Wu, Mingyan},
   title={Schauder estimates for nonlocal equations with singular L\'{e}vy
   Measures},
   journal={Potential Anal.},
   volume={61},
   date={2024},
   number={1},
   pages={13--33},
   issn={0926-2601},
   review={\MR{4758470}},
   doi={10.1007/s11118-023-10101-9},
   eprint={2002.09887},
}

\bib{HZ23}{article}{
   author={Hao, Zimo},
   author={Zhang, Xicheng},
   title={SDEs with supercritical distributional drifts},
   journal={Comm. Math. Phys.},
   volume={406},
   date={2025},
   number={10},
   pages={Paper No. 250, 56},
   issn={0010-3616},
   review={\MR{4952103}},
   doi={10.1007/s00220-025-05430-2},
}
 

\bib{Hol22}{article}{
   author={Holland, Teodor},
   title={On the weak rate of convergence for the Euler-Maruyama scheme with
   H\"older drift},
   journal={Stochastic Process. Appl.},
   volume={174},
   date={2024},
   pages={Paper No. 104379, 16},
   issn={0304-4149},
   review={\MR{4744978}},
   doi={10.1016/j.spa.2024.104379},
}


\bib{HLM17}{article}{
   author={Hu, Y.},
   author={L\^{e}, K.},
   author={Mytnik, L.},
   title={Stochastic differential equation for Brox diffusion},
   journal={Stochastic Process. Appl.},
   volume={127},
   date={2017},
   number={7},
   pages={2281--2315},
   issn={0304-4149},
   review={\MR{3652414}},
   doi={10.1016/j.spa.2016.10.010},
}


\bib{IW89}{book}{
      author={Ikeda, Nobuyuki},
      author={Watanabe, Shinzo},
       title={Stochastic differential equations and diffusion processes},
     edition={Second},
      series={North-Holland Mathematical Library},
   publisher={North-Holland Publishing Co., Amsterdam; Kodansha, Ltd., Tokyo},
        date={1989},
      volume={24},       ISBN={0-444-87378-3},     review={\MR{1011252}},
}

\bib{KP20}{article}{
    AUTHOR = {Kremp, Helena}
  author={Perkowski, Nicolas},
     TITLE = {Multidimensional {SDE} with distributional drift and {L}\'{e}vy
              noise},
   JOURNAL = {Bernoulli},
  FJOURNAL = {Bernoulli. Official Journal of the Bernoulli Society for
              Mathematical Statistics and Probability},
    VOLUME = {28},
      YEAR = {2022},
    NUMBER = {3},
     PAGES = {1757--1783},
      ISSN = {1350-7265},
   MRCLASS = {60H10 (60G51)},
  review = {\MR{4411510}},
       DOI = {10.3150/21-bej1394},
}

  
\bib{LZ22}{article}{
   author={Ling, Chengcheng},
   author={Zhao, Guohuan},
   title={Nonlocal elliptic equation in H\"{o}lder space and the martingale
   problem},
   journal={J. Differential Equations},
   volume={314},
   date={2022},
   pages={653--699},
   issn={0022-0396},
   review={\MR{4369182}},
   doi={10.1016/j.jde.2022.01.025},
}
  
\bib{MP91}{article}{
   author={Mikulevi$\check{\rm c}$ius, Remigius},
   author={Platen, Eckhard},
   title={Rate of convergence of the Euler approximation for diffusion
   processes},
   journal={Math. Nachr.},
   volume={151},
   date={1991},
   pages={233--239},
   issn={0025-584X},
   review={\MR{1121206}},
   doi={10.1002/mana.19911510114},
}
   
\bib{Sa99}{book}{
      author={Sato, Ken-iti},
       title={L\'{e}vy processes and infinitely divisible distributions},
      series={Cambridge Studies in Advanced Mathematics},
   publisher={Cambridge University Press, Cambridge},
        date={1999},
      volume={68},     ISBN={0-521-55302-4},
        note={Translated from the 1990 Japanese original, Revised by the
  author},     review={\MR{1739520}},
}
 


\bib{SH24}{article}{
   author={Song, Ke},
   author={Hao, Zimo},
   title={convergence rates of the Euler-Maruyama scheme to density dependent
   SDEs driven by $\alpha$-stable additive noise},
   journal={Proc. Amer. Math. Soc.},
   volume={153},
   date={2025},
   number={6},
   pages={2591--2607},
   issn={0002-9939},
   review={\MR{4892630}},
   doi={10.1090/proc/17169},
}

\bib{TT90}{article}{
   author={Talay, Denis},
   author={Tubaro, Luciano},
   title={Expansion of the global error for numerical schemes solving
   stochastic differential equations},
   journal={Stochastic Anal. Appl.},
   volume={8},
   date={1990},
   number={4},
   pages={483--509 (1991)},
   issn={0736-2994},
   review={\MR{1091544}},   doi={10.1080/07362999008809220},
}

\bib{Tr92}{book}{
      author={Triebel, Hans},
       title={Theory of function spaces. {II}},
      series={Monographs in Mathematics},
   publisher={Birkh\"{a}user Verlag, Basel},
        date={1992},
      volume={84},     ISBN={3-7643-2639-5},
         url={https://doi.org/10.1007/978-3-0346-0419-2},    review={\MR{1163193}},
}
 
\bib{We19}{article}{
   author={Webb, J. R. L.},
   title={Weakly singular Gronwall inequalities and applications to
   fractional differential equations},
   journal={J. Math. Anal. Appl.},
   volume={471},
   date={2019},
   number={1-2},
   pages={692--711},
   issn={0022-247X},   review={\MR{3906348}},
   doi={10.1016/j.jmaa.2018.11.004},
} 
 
 \bib{Zh10}{article}{
      author={Zhang, Xicheng},
       title={Stochastic {V}olterra equations in {B}anach spaces and stochastic
  partial differential equation},
        date={2010},       ISSN={0022-1236},
     journal={J. Funct. Anal.},
      volume={258},
      number={4},       pages={1361\ndash 1425},
         url={https://doi.org/10.1016/j.jfa.2009.11.006},     review={\MR{2565842}},
} 

\end{biblist}



\end{document}